\documentclass[reqno,centertags]{amsart}
\usepackage{amsmath}
\usepackage{amssymb}
\usepackage{amsfonts}

\newtheorem{theorem}{Theorem}
\theoremstyle{plain}

\newtheorem{corollary}{Corollary}

\newtheorem{lemma}{Lemma}

\newtheorem{problem}{Problem}

\numberwithin{equation}{section}

\input{tcilatex}

\begin{document}
\title[]{Tubes of coordinate finite type Gauss map in the Euclidean 3-space}
\author{Hassan Al-Zoubi}
\address{Department of Mathematics, Al-Zaytoonah University of Jordan, P.O.
Box 130, Amman, Jordan 11733}
\email{dr.hassanz@zuj.edu.jo}
\author{Hamza Alzaareer}
\address{Department of Mathematics, Al-Zaytoonah University of Jordan, P.O.
Box 130, Amman, Jordan 11733}
\email{h.alzaareer@zuj.edu.jo}
\author{Tareq Hamadneh}
\address{Department of Mathematics, Al-Zaytoonah University of Jordan, P.O.
Box 130, Amman, Jordan 11733}
\email{t.hamadneh@zuj.edu.jo}
\author{Mohammad Al Rawajbeh}
\address{Department of Computer Science, Al-Zaytoonah University of Jordan, P.O.
Box 130, Amman, Jordan 11733}
\email{m.rawajbeh@zuj.edu.jo}
\date{}

\begin{abstract}
In this paper, we consider tubes in the Euclidean 3-space whose Gauss map $\boldsymbol{n}$ is of coordinate finite $I$-type, i.e., the  position vector $\boldsymbol{n}$ satisfies the relation $\Delta ^{I}\boldsymbol{n}=\varLambda \boldsymbol{n}$, where $\Delta ^{I}$ is the Laplace operator with respect to the first fundamental form $I$ of the surface and $\varLambda$ is a square matrix of order 3. We show that circular cylinders are the only class of surfaces mentioned above of coordinate finite $I$-type Gauss map.
\end{abstract}

\maketitle

\section{Introduction}
The theory of surfaces of finite Chen type regarding to the first fundamental form $I$ was introduced by B.-Y. Chen about four decades ago and it has been a topic of active research by many differential geometers since then. Many results in this area have been collected in \cite{C7}. A surface $S$ is said to be of finite type corresponding to the first fundamental form $I$, or briefly of finite $I$-type, if the position vector $\boldsymbol{x}$ of $S$ can be written as a finite sum of nonconstant eigenvectors of the Laplacian $\triangle ^{I}$, that is,
\begin{equation*}
\boldsymbol{x}=\boldsymbol{c}+\sum_{i=1}^{k}\boldsymbol{x}_{i},\quad
\Delta ^{I}\boldsymbol{x}_{i}=\lambda _{i}\,\boldsymbol{x}_{i},\quad
i=1,\dotsc ,k,
\end{equation*}%
where $\boldsymbol{c}$ is a fixed vector and $\lambda _{1},\lambda_{2},\dotsc ,\lambda _{k}$ are eigenvalues of the operator $\triangle ^{I}$. In particular, if all eigenvalues $\lambda _{1},\lambda _{2},\dotsc ,\lambda_{k}$ are mutually distinct, then $S$ is said to be of finite $I$-type $k$. When $\lambda_{i}=0$ for some $i=1,\dotsc,k$, then $S$ is said to be of finite null $I$-type $k$. Otherwise, $S$ is said to be of infinite type \cite{C3}.

Denote by $H$ the mean curvature of $S$. Then, it is well known that \cite{T1}
\begin{equation}  \label{2}
\Delta ^{I}\boldsymbol{x}=-2H\boldsymbol{n},
\end{equation}
where the map $\boldsymbol{n} : S \rightarrow M^{2}$ which sends each point of $S$ to the unit normal vector to $S$ at the point is called the Gauss map of the surface $S$, where $M^{2}$ is the unit sphere in the Euclidean 3-space $\mathbb{E}^{3}$ centered at the origin. From (\ref{2}) we know the following two facts \cite{T1}
\begin{itemize}
\item $S$ is minimal if and only if all coordinate functions of $\boldsymbol{x}$ are eigenfunctions of $\Delta ^{I}$ with eigenvalue $\lambda=0$.

\item $S$ lies in an ordinary sphere $M^{2}$ if and only if all coordinate functions of $\boldsymbol{x}$ are eigenfunctions of $\Delta ^{I}$ with a fixed nonzero eigenvalue.
\end{itemize}

Results concerning surfaces of finite type in $\mathbb{E}^{3}$ remain very little. In fact, spheres, minimal surfaces and circular cylinders are the only known surfaces of finite $I$-type in $\mathbb{E}^{3}$. Therefore in \cite{C2} B.-Y. Chen posed the following interesting question
\begin{problem}  \label{(p1)}
Determine all surfaces of finite $I$-type in $\mathbb{E}^{3}$.
\end{problem}
Many researchers start to solve this problem by investigating important classes of surfaces in $\mathbb{E}^{3}$. More specifically, ruled surfaces \cite{C5}, quadrics \cite{C6}, tubes \cite{C4}, cyclides of Dupin \cite{D1, D2} and spiral surfaces \cite{B1} are the only known families of surfaces in $\mathbb{E}^{3}$ that have been studied according to its finite type classification. However, for another important families of surfaces, such as translation surfaces, helicoidal surfaces as well as surfaces of revolution, the classification of its finite type still unknown. For a more details, the reader can refer to \cite{C7}.

Later in \cite{G3} O. Garay generalized T. Takahashi's condition studied surfaces in $\mathbb{E}^{3}$ for which all coordinate functions $\left(x_{1}, x_{2}, x_{3}\right)$ of $\boldsymbol{x}$ satisfy $\Delta^{I}\boldsymbol{x}_{i} = \lambda_{i}x_{i}, i = 1,2,3$, not necessarily with the same eigenvalue. Another generalization was studied in \cite{D3} for which surfaces in $\mathbb{E}^{3}$ satisfy the relation $\Delta^{I}\boldsymbol{x}= B\boldsymbol{x} + C$ $(*)$ where $B \in\mathbb{Re}^{3\times3} ;C \in\mathbb{Re}^{3\times1}$. It was shown that a surface $S$ in $\mathbb{E}^{3}$ satisfies $(*)$ if and only if it is an open part of a minimal surface, a sphere, or a circular cylinder. A surface $S$ in $\mathbb{E}^{3}$ whose position vector $\boldsymbol{x}$ satisfies the condition $(*)$ is said to be of coordinate finite type.

Another interesting theme within this context is to study families of surfaces in the Euclidean 3-space for which its Gauss map is of finite type \cite{C9}. Here again, we also have the following question
\begin{problem}  \label{(p2)}
Which surfaces in $\mathbb{E}^{3}$ are of finite type Gauss map.
\end{problem}
Concerning this question, in \cite{B4} it was proved that
\begin{theorem}\label{TH1}
The only ruled surfaces in the $n$-dimensional Euclidean space $(n\geq 3)$ with finite type Gauss map are cylinders over curves of finite type and planes.
\end{theorem}
\begin{theorem}
\label{TH2} The only tubes in $\mathbb{E}^{3}$ of finite type Gauss map are circular cylinders.
\end{theorem}

In \cite{B3} C. Baikoussis and others studied the class of cyclides of Dupin, they proved that the compact, as well as for the noncompact cyclides of Dupin, the Gauss map is of infinite type. Further, in \cite{B1} the problem was solved for the family of spiral surfaces, it was shown that planes are the only spiral surfaces in $\mathbb{E}^{3}$ with finite type Gauss map. Meanwhile, the second problem still unsolved for the quadric surfaces, the surfaces of revolution, the cones, the translation surfaces and the helicoidal surfaces.

Following the ideas of O. Garay in \cite{G3}, one can study all classes of surfaces in $\mathbb{E}^{3}$ mentioned above whose Gauss map satisfies the following condition
\begin{equation}  \label{3}
\Delta ^{I}\boldsymbol{n}=\Lambda \boldsymbol{n},
\end{equation}%
where $\Lambda \in\mathbb{Re}^{3\times3}$. Therefore, an interesting geometric question is raised:
\begin{problem}  \label{p3}
Determine all surfaces in $\mathbb{E}^{3}$ with the Gauss map satisfying (\ref{3}).
\end{problem}

In this respect, the class of surfaces of revolution was studied in \cite{D4}, it was shown that planes, spheres and circular cylinders are the only surfaces of revolution whose Gauss map satisfies (\ref{3}). In \cite{B5} Ch. Baikoussis and L. Verstraelen studied the class of translation surfaces and they proved that the only translation surfaces whose Gauss map satisfies (\ref{3}) are the planes and the circular cylinders. The same authors in \cite{B6} solved this problem for the class of helicoidal surfaces, they proved that planes, spheres and circular cylinders are the only helicoidal surfaces whose Gauss map satisfies condition (\ref{3}). Moreover, in \cite{B2} Ch. Baikoussis and D. Blair showed that circular cylinders and planes are the only ruled surfaces for which its Gauss map satisfies (\ref{3}). In \cite{B3} it was proved that neither for the compact, nor for the noncompact cyclides of Dupin, condition (\ref{3}) can be satisfied. Finally, in \cite{B1} it was proved that planes are the only spiral surfaces in $\mathbb{E}^{3}$ whose Gauss map satisfies (\ref{3}). Here again the classification of quadric surfaces, tubes and cones in the sense of condition (\ref{3}) still unsolved.

Sixteen years ago, in \cite{S1} S. Stamatakis and H. Al-Zoubi introduced the notion of surfaces of finite type in $\mathbb{E}^{3}$ with respect to the second or third fundamental form of a surface and it has become a useful tool from that time on as one can see in the recent literature in this field. In this respect, it is worthwhile to study all the problems mentioned above concerning the second or third fundamental form of a surface. Many results in this area can be found in \cite{A1, A2, A3, A4, A5, A6, A7, K1, S0, S2}.

In this paper we will focus on problem (\ref{p3}) by studying an important class of surfaces, namely, tubes in $\mathbb{E}^{3}$. Our main theorem is:
\begin{theorem}
\label{TH3} Circular cylinders are the only tubes in the 3-dimensional Euclidean space whose Gauss map $\boldsymbol{n}$ satisfying (\ref{3}).
\end{theorem}




\section{Tubes in $\mathbb{E}^{3}$}

Let $C: \boldsymbol{\alpha} =\boldsymbol{\alpha}(\emph{t})$, $\mathit{t}\epsilon (a,b)$ be a regular unit speed curve of finite length which is topologically imbedded in $\mathbb{E}^{3}$. The total space $N_{\boldsymbol{\alpha}}$ of the normal bundle of $\boldsymbol{\alpha}(a, b)$ in $\mathbb{E}^{3}$ is naturally diffeomorphic to the direct product $(a,b)\times \mathbb{E}^{2}$ via the translation along $\boldsymbol{\alpha}$ with respect to the induced normal connection. For a sufficiently small $r>0$ the tube of radius $r$ about the curve $\boldsymbol{\alpha}$ is the set
\begin{equation}
T_{r}( \boldsymbol{\alpha})=\{exp_{\boldsymbol{\alpha}(t)}\boldsymbol{u}\mid \boldsymbol{u}\in N_{\boldsymbol{\alpha}} , \ \ \parallel \boldsymbol{u}%
\parallel =r,\ \ t\in(a,b)\}.  \notag
\end{equation}
Assume that ${\mathbf{t}, \mathbf{h}, \mathbf{b}}$ is the Frenet frame, and $\kappa$ is the curvature of the unit speed curve $\boldsymbol{\alpha} =
\boldsymbol{\alpha}(\emph{t})$. For a small real number $r$ satisfies $0 < r< min\frac{1}{|\kappa|}$, the tube $T_{r}( \boldsymbol{\alpha})$ is a smooth
surface in $\mathbb{E}^{3}$ \cite{R1}. Then a parametric representation of the tube $T_{r}( \boldsymbol{\alpha})$ is given by

\begin{equation*}  \label{eq6}
\digamma:\boldsymbol{x}(t,\varphi)= \boldsymbol{\alpha}+ r \cos\varphi\mathbf{h}+ r \sin\varphi\mathbf{b}.
\end{equation*}

The first fundamental form of $\digamma$ is
\begin{align*}
I &= \big(\delta^{2}+r^{2}\tau^{2}\big)dt^{2} + 2r^{2}\tau dtd\varphi+r^{2}d\varphi^{2},
\end{align*}
where $\delta: = (1-r\kappa\cos\varphi)$ and $\tau$ is the torsian of the curve $\boldsymbol{\alpha}$. The Gauss map $\boldsymbol{n}$ of $\digamma$ is given by

\begin{equation}  \label{eq71}
\boldsymbol{n} =-(\cos\varphi \mathbf{h}+\sin\varphi\mathbf{b}).
\end{equation}
The Beltrami-Laplace operator corresponding to the first fundamental form of $\digamma$ can be expressed as follows \cite{C4}

\begin{eqnarray}  \label{eq8}
\Delta ^{I} &=&-\frac{1}{\delta^{3}}\Bigg[\delta\frac{\partial^{2}}{\partial t^{2}}-2\tau \delta\frac{\partial ^{2}}{\partial t\partial \varphi}%
+\frac{\delta}{r^{2}}(r^{2}\tau ^{2}+\delta^{2})\frac{\partial^{2}}{\partial \varphi^{2}}\\ \notag
&&+r\beta\frac{\partial}{\partial t}-\frac{\kappa \delta^{2}\sin\varphi}{r}\frac{\partial }{\partial \varphi }\Bigg] ,
\end{eqnarray}
where $\beta: =\kappa \acute{} \cos \varphi +\kappa \tau \sin \varphi $ and $\acute{}:=\frac{d}{dt}$.

Before we start proving of the main result, we mention and prove the following two special cases of tubular surfaces for later use.

\subsection{Anchor rings}

A tube in the Euclidean 3-space is called an anchor ring if the curve $C$ is a plane circle (or is an open portion of a plane circle). In this case, the torsian $\tau$ of $\boldsymbol{\alpha}$ vanishes identically, and the curvature $\kappa$ of $\boldsymbol{\alpha}$ is a nonzero constant. Then the position vector $\boldsymbol{x}$ of the anchor ring can be expressed as

\begin{equation*}  \label{eq9}
\digamma:\boldsymbol{x}(t,\varphi)= \{\gamma \cos\varphi,\gamma \sin\varphi, r\sin t \},
\end{equation*}
where $\gamma :=a+r\cos t, a > r$ and $a \epsilon \mathbb{R}$. The first fundamental form is

\begin{equation*}
I=r^{2}dt^{2}+\gamma^{2}d\varphi ^{2}.
\end{equation*}

Hence, the Beltrami-Laplace operator is given by  \cite{C4}

\begin{equation}  \label{eq10}
\Delta^{I}= -\frac{1}{\gamma^{^{2}}}\frac{\partial^{2}}{\partial \varphi^{2}}+\frac{\sin t}{r\gamma}\frac{\partial}{\partial t} - \frac{1}{r^{^{2}}}\frac{\partial^{2}}{\partial t^{2}}.
\end{equation}

Denoting by $\boldsymbol{n}$ the Gauss map of $\digamma$, then we have
\begin{equation*}
\boldsymbol{n}=\{-\cos t \cos \varphi , -\cos t \sin\varphi , -\sin t\}.
\end{equation*}

Let $(n_{1}, n_{2}, n_{3})$ be the coordinate functions of $\boldsymbol{n}$. By virtue of (\ref{eq10}) one can find

\begin{equation}  \label{eq11}
\Delta^{I}n_{1} = \Bigg[\frac{\sin^{2}t}{r\gamma} -\Big(\frac{1}{r^{2}}+\frac{1}{\gamma^{2}}\Big)\cos t\Bigg]\cos\varphi,
\end{equation}

\begin{equation}  \label{eq12}
\Delta^{I}n_{2} = \Bigg[\frac{\sin^{2}t}{r\gamma} -\Big(\frac{1}{r^{2}}+\frac{1}{\gamma^{2}}\Big)\cos t\Bigg]\sin\varphi,
\end{equation}

\begin{equation}  \label{eq13}
\Delta^{I}n_{3} = -\frac{\sin t}{r}\Bigg[\frac{\cos t}{\gamma} +\frac{1}{r}\Bigg].
\end{equation}

We denote the entries of the matrix $\Lambda $ by $\lambda _{ij}$ for $i,j=1,2,3$. On account of (\ref{3}) and from (\ref{eq11}), (\ref{eq12}) and (\ref{eq13}) we get

\begin{equation}  \label{eq14}
\Bigg[\frac{\sin^{2}t}{r\gamma} -\Big(\frac{1}{r^{2}}+\frac{1}{\gamma^{2}}\Big)\cos t\Bigg]\cos\varphi=-\lambda _{11}\cos t \cos \varphi-\lambda _{12}\cos t \sin\varphi-\lambda _{13}\sin t,
\end{equation}

\begin{equation}  \label{eq15}
\Bigg[\frac{\sin^{2}t}{r\gamma} -\Big(\frac{1}{r^{2}}+\frac{1}{\gamma^{2}}\Big)\cos t\Bigg]\sin\varphi=-\lambda _{21}\cos t \cos \varphi-\lambda _{22}\cos t \sin\varphi-\lambda _{23}\sin t,
\end{equation}%

\begin{equation}  \label{eq16}
-\frac{\sin t}{r}\Bigg[\frac{\cos t}{\gamma} +\frac{1}{r}\Bigg]=-\lambda _{31}\cos t \cos \varphi-\lambda _{32}\cos t \sin\varphi-\lambda _{33}\sin t.
\end{equation}

Differentiating (\ref{eq14}) and (\ref{eq15}) twice with respect to $\varphi$, it is easy to see that $\lambda _{13}=\lambda _{23} =0$. Similarly, if we take the derivative of (\ref{eq16}) with respect to $\varphi$, then we have $\lambda _{31}=\lambda _{32} =0$. So relations (\ref{eq14}), (\ref{eq15}) and (\ref{eq16}) reduce to

\begin{equation*}  \label{eq17}
\Bigg[\frac{\sin^{2}t}{r\gamma} -\Big(\frac{1}{r^{2}}+\frac{1}{\gamma^{2}}-\lambda _{11}\Big)\cos t\Bigg]\cos\varphi +\lambda _{12}\cos t \sin\varphi =0,
\end{equation*}

\begin{equation*}  \label{eq18}
\Bigg[\frac{\sin^{2}t}{r\gamma} -\Big(\frac{1}{r^{2}}+\frac{1}{\gamma^{2}}-\lambda _{22}\Big)\cos t\Bigg]\sin\varphi +\lambda _{21}\cos t \cos \varphi = 0,
\end{equation*}%

\begin{equation*}  \label{eq19}
-\frac{\sin t}{r}\Bigg[\frac{\cos t}{\gamma} +\frac{1}{r}\Bigg]+\lambda _{33}\sin t=0.
\end{equation*}

Since $\sin\varphi$ and $\cos \varphi$ are linearly independent functions, one finds that $\lambda _{21}=\lambda _{12} =0$. Hence, we have the following three equations

\begin{equation*}
[r\gamma\sin^{2}t -(r^{2}+\gamma^{2}-\lambda _{11}r^{2}\gamma^{2})\cos t]\cos\varphi =0,
\end{equation*}

\begin{equation*}
[r\gamma\sin^{2}t -(r^{2}+\gamma^{2}-\lambda _{22}r^{2}\gamma^{2})\cos t]\sin\varphi = 0,
\end{equation*}

\begin{equation*}
-r\sin t\cos t -\gamma\sin t+\lambda _{33}r^{2}\gamma\sin t=0.
\end{equation*}

From the above three equation, we obtain
\begin{equation*}  \label{eq20}
r\gamma\sin^{2}t -(r^{2}+\gamma^{2}-\lambda _{11}r^{2}\gamma^{2})\cos t =0,
\end{equation*}

\begin{equation*}  \label{eq21}
r\gamma\sin^{2}t -(r^{2}+\gamma^{2}-\lambda _{22}r^{2}\gamma^{2})\cos t = 0,
\end{equation*}%

\begin{equation*}  \label{eq22}
-r\cos t -\gamma+\lambda _{33}r^{2}\gamma=0.
\end{equation*}

This is impossible since we have $\lambda _{11}, \lambda _{22}$ and $\lambda _{33}$ depend on the parameter $t$ and they are not constants.
Consequently, we have the following

\begin{corollary}
\label{C1.1} For the anchor ring, there exists no matrix $\Lambda \in\mathbb{Re}^{3\times3}$ such that the condition (\ref{3}) is satisfied.
\end{corollary}

\subsection{Circular cylinder}
A tube in the Euclidean 3-space is called a circular cylinder when the curve $\boldsymbol{\alpha}$ lies in a line. In this case, the curvature $\kappa$ of $\boldsymbol{\alpha}$ vanishes identically. The position vector $\boldsymbol{x}$ of the circular cylinder can be expressed as

\begin{equation}  \label{eq23}
\digamma:\boldsymbol{x}(t,\varphi)= \{r\cos\varphi ,r\sin\varphi ,t\},   \ \ \ r>0.
\end{equation}

The Gauss map of $\digamma$ is given by
\begin{equation}  \label{eq24}
\boldsymbol{n}(t,\varphi)= -\{\cos\varphi,\sin\varphi , 0\}.
\end{equation}

Since $\kappa =0$, we have $\delta=1$, and so relation (\ref{eq8}) takes the following form
\begin{eqnarray}  \label{eq25}
\Delta ^{I} &=&-\frac{\partial^{2}}{\partial t^{2}}-\frac{1}{r^{2}}\frac{\partial^{2}}{\partial \varphi^{2}} ,
\end{eqnarray}

By using (\ref{eq24}) and (\ref{eq25}) we find

\begin{eqnarray*}  \label{eq26}
\Delta ^{I}\boldsymbol{n} =\frac{1}{r^{2}}\boldsymbol{n},
\end{eqnarray*}
and the following lemma holds
\begin{lemma}
\label{l1.1} The Gauss map of a circular cylinder in the Euclidean 3-space with parametric representation (\ref{eq23}) satisfies (\ref{3}). The corresponding matrix is
\begin{equation*}
\Lambda =\left[
\begin{array}{ccc}
\frac{1}{r^{2}} & 0 & \lambda_{13} \\
0 & \frac{1}{r^{2}} & \lambda_{23} \\
0 & 0& \lambda_{33}
\end{array}%
\right],
\end{equation*}
\end{lemma}
\noindent where $\lambda_{i3}, i=1,2,3$, are arbitrary constants \cite{B2}.

\section{Proof of the main theorem}

Applying relation (\ref{eq8}) on the position vector $\boldsymbol{n}$ of (\ref{eq71}) gives

\begin{eqnarray*}  \label{eq27}
\Delta ^{I}\boldsymbol{n} &=& -\frac{\beta}{\delta^{3}}\mathbf{t}-\frac{1}{\delta^{2}}\big(\kappa^{2}\cos\varphi -2\tau^{2}\cos\varphi +\frac{\delta^{2}}{r^{2}}\cos\varphi +\frac{\kappa\delta}{r}\sin^{2}\varphi \big)\mathbf{h}\\  \notag
&&+\frac{1}{r^{2}\delta}\sin\varphi(2r\kappa\cos\varphi-1)\mathbf{b}.
\end{eqnarray*}

Let $\{t_{1}, t_{2}, t_{3}\}, \{h_{1}, h_{2}, h_{3}\}$ and $\{b_{1}, b_{2}, b_{3}\}$ be the components of $\mathbf{t},\mathbf{h}$ and $\mathbf{b}$ respectively. On account of (\ref{3}) we get

\begin{eqnarray*}  \label{eq28}
&&-\frac{\beta}{\delta^{3}}t_{i}-\frac{1}{\delta^{2}}\big(\kappa^{2}\cos\varphi -2\tau^{2}\cos\varphi +\frac{\delta^{2}}{r^{2}}\cos\varphi +\frac{\kappa\delta}{r}\sin^{2}\varphi \big)h_{i}\\ \notag
&&+\frac{1}{r^{2}\delta}\sin\varphi(2r\kappa\cos\varphi-1)b_{i}=\\  \notag
&& -\lambda_{i1}(\cos\varphi h_{1}+\sin\varphi b_{1}) - \lambda_{i2}(\cos\varphi h_{2}+\sin\varphi b_{2})\\ \notag
&&-\lambda_{i3}(\cos\varphi h_{3}+\sin\varphi b_{3}), \ \ i= 1, 2, 3, \notag
\end{eqnarray*}
which can be rewritten as
\begin{eqnarray}  \label{eq29}
&&-\beta t_{i}-\delta \big(\kappa^{2}\cos\varphi -2\tau^{2}\cos\varphi +\frac{\delta^{2}}{r^{2}}\cos\varphi +\frac{\kappa\delta}{r}\sin^{2}\varphi \big)h_{i}\\ \notag
&&+\frac{\delta^{2}}{r^{2}}\sin\varphi(2r\kappa\cos\varphi-1)b_{i}=\\  \notag
&& -\lambda_{i1}\delta^{3}(\cos\varphi h_{1}+\sin\varphi b_{1}) - \lambda_{i2}\delta^{3}(\cos\varphi h_{2}+\sin\varphi b_{2})\\ \notag
&&-\lambda_{i3}\delta^{3}(\cos\varphi h_{3}+\sin\varphi b_{3}). \notag
\end{eqnarray}

For $i = 1, 2, 3$, the left hand side of (\ref{eq29}) is a polynomial in $\cos\varphi$ and $\sin\varphi$ with functions in $t$ as coefficients. This implies that the coefficients of the powers of $\cos\varphi$ and $\sin\varphi$ must be zeros, so we obtain, for $i = 1, 2, 3$, the following equations

\begin{eqnarray}  \label{eq30}
&&-\beta t_{i}-\delta \big(\kappa^{2}\cos\varphi -2\tau^{2}\cos\varphi +\frac{\delta^{2}}{r^{2}}\cos\varphi +\frac{\kappa\delta}{r}\sin^{2}\varphi \big)h_{i}\\ \notag
&&+\frac{\delta^{2}}{r^{2}}\sin\varphi(2r\kappa\cos\varphi-1)b_{i}= 0,  \notag
\end{eqnarray}
and
\begin{eqnarray*}  \label{eq31}
&& \lambda_{i1}\delta^{3}(\cos\varphi h_{1}+\sin\varphi b_{1}) + \lambda_{i2}\delta^{3}(\cos\varphi h_{2}+\sin\varphi b_{2})\\ \notag
&&+\lambda_{i3}\delta^{3}(\cos\varphi h_{3}+\sin\varphi b_{3})=0. \notag
\end{eqnarray*}

From (\ref{eq30}), we conclude
\begin{eqnarray*}  \label{eq32}
&&-\beta \mathbf{t}-\delta \big(\kappa^{2}\cos\varphi -2\tau^{2}\cos\varphi +\frac{\delta^{2}}{r^{2}}\cos\varphi +\frac{\kappa\delta}{r}\sin^{2}\varphi \big)\mathbf{h}\\ \notag
&&+\frac{\delta^{2}}{r^{2}}\sin\varphi(2r\kappa\cos\varphi-1)\mathbf{b}= 0,
\end{eqnarray*}
which implies that $\beta=0$. Therefore $\kappa'= 0$ and $\kappa\tau = 0$, so we have the following two cases:

Case 1. $\kappa$ = const. $\neq0$ and $\tau = 0$. Thus the curve $C$ is a plane circle, and so $\digamma$ is anchor ring. Hence,
according to Corollary \ref{C1.1}, there exists no matrix $\Lambda \in\mathbb{Re}^{3\times3}$ such that the condition (\ref{3}) is satisfied.

Case 2. $\kappa=0$. Hence the curve $C$ is a line. Therefore the torsian $\tau$ of the curve $C$ vanishes identically, and so according to Lemma \ref{l1.1}, $\digamma$ is circular cylindar $\Box$  




\end{document}